\documentclass[12pt]{article}
\usepackage{cite,amsfonts,amsmath,subeqn}
\usepackage{epic}
\def\be{\begin{equation}}
\def\eqn#1{\be\label{#1}}

\def\bea{\begin{eqnarray}}

\def\eea{\end{eqnarray}}

\newcommand{\eqna}[1]{\begin{subequations} \label{#1}
\begin{eqnarray}}
\def\eena{\end{eqnarray}
\end{subequations}}

\def\id{{\bf 1}}

\def\nl{\hfil\break}
\def\nt{\noindent}

\def\nn{\nonumber}
\def\eps{{\epsilon}}

\newtheorem{definition}{Definition}[section]
\newtheorem{example}[definition]{Example}

\newtheorem{proposition}[definition]{Proposition}

\def\diag{\mathop{\rm diag}\nolimits}
\newcommand{\qmbox}[1]{{\qquad\mbox{#1}\quad}}

\def\d{{\delta}}
  \def\cC{{\cal C}}
  
 \def\cH{{\cal H}} 
 \def\cK{{\cal K}} 
  
  \def\cR{{\cal R}}
 \def\cT{{\cal T}}

  \def\tL{\tilde L}

\def\hR{\hat R}
\def\up{{\uparrow}}
\def\dw{{\downarrow}}
\def\state#1{\left| #1 \right>}

\normalbaselines
\begin{document}

 \begin{center}

 \textsf{\LARGE Exotic Bialgebra $S03$: 
Representations, Baxterisation and Applications}

 \vspace{10mm}

 {\large D.~Arnaudon$^{a,}$\footnote{daniel.arnaudon@lapp.in2p3.fr},
 ~A.~Chakrabarti$^{b,}$\footnote{chakra@cpht.polytechnique.fr},\\[2mm]
 V.K.~Dobrev$^{c,}$\footnote{dobrev@inrne.bas.bg}
 ~and~ S.G.~Mihov$^{c,}$\footnote{smikhov@inrne.bas.bg}}

 \vspace{5mm}

 \emph{$^a$ Laboratoire d'Annecy-le-Vieux de Physique Th{\'e}orique LAPTH}
 \\
 \emph{CNRS, UMR 5108, associ{\'e}e {\`a} l'Universit{\'e} de Savoie}
 \\
 \emph{LAPTH, BP 110, F-74941 Annecy-le-Vieux Cedex, France}
 \\
 \vspace{3mm}
 \emph{$^b$ Centre de Physique Th{\'e}orique, CNRS UMR 7644}
 \\
 \emph{Ecole Polytechnique, 91128 Palaiseau Cedex, France.}
 \\
 \vspace{3mm}
 \emph{$^c$ Institute of Nuclear Research and Nuclear Energy}
 \\
 \emph{Bulgarian Academy of Sciences}
 \\
 \emph{72 Tsarigradsko Chaussee, 1784 Sofia, Bulgaria}
 \\
 \vspace{3mm}
\end{center}

\begin{abstract}
The exotic bialgebra  $S03$, defined by a solution of the
Yang-Baxter  equation, which is not a deformation of the trivial,
is considered. Its FRT dual algebra $s03_F$ is  studied.
The Baxterisation of the dual algebra is given in two different
parametrisations. The finite-dimensional representations of $s03_F$
are considered. Diagonalisations of the braid matrices are used to yield
remarkable insights concerning representations of the $L$-algebra
and to formulate the fusion of finite-dimensional representations.
Possible applications are considered, in particular, an exotic
eight-vertex  model and  an integrable spin-chain model.
\end{abstract}

\vfill \centerline{\textit{Dedicated to our friend Daniel Arnaudon}}

\vfill

\vfill

\noindent
MSC: 20G42, 17B37, 81R50 \hfill
{LAPTH-1137/06}\\
PACS: 02.20.Uw, 02.10.Hh,  02.30.Ik\hfill {math.QA/0601708}
\rightline{January 2006}
\rightline{revised March 2006}

\vfill\eject

\section{Introduction}
\label{sect:intro}
\setcounter{equation}{0}

For several years \cite{ACDM1,ACDM2,ACDM3}
our collaboration studied the   algebraic structures coming from
4x4 $R$-matrices (solutions of the Yang--Baxter equation)
that are not deformations of classical ones (i.e., the
identity up to signs). According to
the classification of Hietarinta \cite{Hietarinta}
there are five such 4x4 $R$-matrices
that are invertible.
These matrices were obtained first by Hlavat\'y
\cite{Hlavaty} without classification claims.
In the present paper we consider in more detail one of these cases which seems
most interesting, namely, the matrix bialgebra $S03$ and its FRT \cite{FRT} dual $s03_F\,$.

The paper is organised as follows.
In Section 2 we introduce the matrix bialgebra $S03$, its FRT dual $s03_F\,$,
and an affinisation for the latter.
We give also a basis of $s03_F$ suitable to define the class of representations
when $s03_F$ acts on itself.
In Section 3 we give alternative parametrisation of the baxterised $R$ and $L$ matrices.
This allows to introduce a diagonalisation of the braid matrix which gives
remarkable insights concerning representations of the L-algebra.
In Section 4 we study the finite-dimensional representations of $s03_F\,$.
In Section 5 we use the diagonalisation of the permuted $R$-matrix in order to
formulate the fusion of certain representations.
In Section 6 we consider some of the possible applications:
an exotic eight-vertex  model and  an integrable spin-chain model
are discussed.

\section{FRT Duality}
\label{sect:prelim}
\setcounter{equation}{0}

\subsection{Preliminaries}

Our starting point is the following $4\times 4$ $R$-matrix:
\begin{equation}
 \label{eq:R_S03}
 R = R_{S03} = \frac1{\sqrt2}
 \left(
 \begin{array}{cccc}
 1 & 0 & 0 & 1 \cr
 0 & 1 & 1 & 0 \cr
 0 & 1 & -1 & 0 \cr
 -1 & 0 & 0 & 1
 \end{array}
 \right) \ .
\end{equation}
This $R$-matrix appears in the classification of \cite{Hietarinta}
which gives all (up to equivalence) $4\times 4$ matrix
solutions of the Yang-Baxter equation. Obviously, \eqref{eq:R_S03}
is not a deformation of the identity.\footnote{Higher
dimensional ($N^2 \times N^2$ matrices for all $N>2$) exotic braid matrices (which are
not  deformations of some  "classical limits") have been presented and studied in
\cite{AC0401,Chakra3}.}

In this subsection we introduce various quantities that we need later.
First two standard matrices $R^\pm$  defined by:
\begin{equation}
 \label{eq:Rpm}
 R^+ \equiv PRP  = \frac1{\sqrt2}
 \left(
 \begin{array}{cccc}
 1 & 0 & 0 & 1 \cr
 0 & -1 & 1 & 0 \cr
 0 & 1 & 1 & 0 \cr
 -1 & 0 & 0 & 1
 \end{array}
 \right) \;,\qquad
 R^- \equiv R^{-1} = \frac1{\sqrt2}
 \left(
 \begin{array}{cccc}
 1 & 0 & 0 & -1 \cr
 0 & 1 & 1 & 0 \cr
 0 & 1 & -1 & 0 \cr
 1 & 0 & 0 & 1
 \end{array}
 \right) \ ,
\end{equation}
where ~$P$~ is the permutation matrix: \eqn{perm} P\ \equiv\
\left(
 \begin{array}{cccc}
 1 & 0 & 0 & 0 \cr
 0 & 0 & 1 & 0 \cr
 0 & 1 & 0 & 0 \cr
 0 & 0 & 0 & 1 \cr
 \end{array}
 \right)
\end{equation}

Then we introduce the Baxterised $R$-matrix:
\begin{eqnarray}
 \label{eq:R_S03x}
 R(x) = R_{S03}(x)
 &=&
 x^{-1/2} R + x^{1/2} R_{21}^{-1} =
 \nonumber\\[3mm]
 &=&
 \frac1{\sqrt{2x}}
 \left(
 \begin{array}{cccc}
 x+1 & 0 & 0 & 1-x \cr
 0 & 1-x & x+1 & 0 \cr
 0 & x+1 & x-1 & 0 \cr
 x-1 & 0 & 0 & x+1
 \end{array}
 \right) \ .
\end{eqnarray}
It satisfies the spectral parameter dependent Yang--Baxter equation
\begin{equation}
 \label{eq:pybe}
R_{12}(x)R_{13}(xy)R_{23}(y) =
R_{23}(y)R_{13}(xy)R_{12}(x)
\end{equation}

Finally, we define the braid matrix
~$\hat{R} \equiv PR$, and its Baxterisation $\hat{R}(x)$~:
\begin{eqnarray}
 \label{eq:hR}
 \hR &\equiv& P_{12} R
 \\
 \label{eq:hRx}
 \hR(x) &\equiv& P_{12} R(x) =
 x^{-1/2} \hR + x^{1/2} \hR^{-1} =
 \nonumber\\[3mm]
&=& \frac{1}{\sqrt{2x}}
  \left(
    \begin{array}{cccc}
      x+1 & 0 & 0 & 1-x \cr
      0 & x+1 & x-1 & 0 \cr
      0 & 1-x & x+1 & 0 \cr
      x-1 & 0 & 0 & x+1
    \end{array}
  \right)
\end{eqnarray}
Note that  $\hat{R}^{-1} = R^{-}P$.

We also record two identities involving $\hR$~:
\begin{equation}
 \label{eq:hRprop}
 \hR + \hR^{-1} = \sqrt2 I \qquad\qquad
 \hR^2 + \hR^{-2} = 0
\end{equation}

\subsection{The bialgebra S03}

Here we recall  the matrix bialgebra ~$S03$~ which we obtained in
\cite{ACDM2} by applying the RTT relations of \cite{FRT}:
\eqn{rtt} R\ T_1\ T_2 \ \ =\ \ T_2\ T_1\ R \  , \qquad T = \left(
\begin{array}{cc} a & b \cr c & d \end{array} \right)
\end{equation} where \ $T_1 \ =\ T\, \otimes\, \id_2$\ , \ $T_2 \
=\ \id_2 \, \otimes\, T$, for the case when $\ R\ =\ R_{S0,3}\ $.
The relations which follow from (\ref{rtt}) and (\ref{eq:R_S03})
are:
\begin{alignat}{2}
 \label{eq:S03rel}
 & b^2 + c^2 = 0\ , \qquad&&
 a^2 - d^2 = 0\ , \nn\\
 & cd = ba\ , &&
 dc = -ab\ , \nn\\
 & bd = ca\ , &&
 db = -ac\ , \nn\\
 & da = ad\ , &&
 cb = -bc\ .
\end{alignat}

\subsection{The FRT dual s03$_F$}
\label{sect:RLL}

The FRT dual ~$s03_F$~of ~$S03$~ is given in terms of
$L^\pm$ which are matrices of operators $L^\pm_{ij}$ ($i,j=1,2$)
satisfying the so-called RLL relations \cite{FRT}:
\begin{eqnarray}
 \label{eq:RLL}
 R^+ L^+_1 L^+_2 = L^+_2 L^+_1 R^+ \nn\\
 R^+ L^-_1 L^-_2 = L^-_2 L^-_1 R^+ \nn\\
 R^+ L^+_1 L^-_2 = L^-_2 L^+_1 R^+
\end{eqnarray}
with $L_1\equiv L\otimes 1$, $L_2\equiv 1\otimes L$.

Encoding $L^+$ and $L^-$ in $L(x) = x^{-1/2} L^+ + x^{1/2} L^-$, the
equations (\ref{eq:RLL}) are equivalent to
\begin{equation}
 \label{eq:pybel}
R_{12}(x/y) L_{1}(x) L_{2}(y) =
L_{2}(y) L_{1}(x) R_{12}(x/y)
\end{equation}

Explicitly, these RLL relations read
\begin{alignat}{2}
 \label{eq:RLLij}
 & (L^{\pm}_{11})^2 = (L^{\pm}_{22})^2
 &\qquad\qquad
 & [L^{\pm}_{11}, L^{\pm}_{22} ] = 0 \nonumber\\
 & (L^{\pm}_{12})^2 = - (L^{\pm}_{21})^2
 &
 & [L^{\pm}_{12}, L^{\pm}_{21} ]_+ = 0 \nonumber\\
 & L^{\pm}_{11} L^{\pm}_{12} = L^{\pm}_{22} L^{\pm}_{21}
 &
 & L^{\pm}_{11} L^{\pm}_{21} = L^{\pm}_{22} L^{\pm}_{12} \nonumber\\
 & L^{\pm}_{12} L^{\pm}_{11} = - L^{\pm}_{21} L^{\pm}_{22}
 &
 & L^{\pm}_{12} L^{\pm}_{22} = - L^{\pm}_{21} L^{\pm}_{11}
\end{alignat}
and for the $RL^+ L^-$ ones
\begin{equation}
 \label{eq:RL+L-}
 L^+_{ij} L^-_{kl} - L^-_{ij} L^+_{kl}
 + \theta_i L^+_{\bar ij} L^-_{\bar kl} + \theta_j L^-_{i\bar j}
 L^+_{k\bar l} = 0
\end{equation}
with $\bar n \equiv 3-n$, $\theta_1=1$, $\theta_2=-1$.
The $RL^+L^+$ relations are to be compared with \eqref{eq:S03rel}.

Introducing
\begin{alignat}{2}
 &\tL^\pm_{11} = L^\pm_{11} + L^\pm_{22}
 &\qquad\qquad
 &\tL^\pm_{22} = L^\pm_{11} - L^\pm_{22}
 \nonumber\\
 &\tL^\pm_{12} = L^\pm_{12} + L^\pm_{21}
 &
 &\tL^\pm_{21} = L^\pm_{12} - L^\pm_{21}
\end{alignat}
the relations (\ref{eq:RLLij}) become
\begin{alignat}{2}
 \label{eq:RtLtL1}
 &\tL^\pm_{11} \,\tL^\pm_{22} = 0
 &\qquad\qquad\qquad
 &\tL^\pm_{22} \,\tL^\pm_{11} =0
 \nonumber\\
 &(\tL^\pm_{12}) ^2 = 0
 &
 &(\tL^\pm_{21}) ^2 = 0
 \nonumber\\
 &\tL^\pm_{11} \,\tL^\pm_{21} = 0
 &
 &\tL^\pm_{12} \,\tL^\pm_{11} =0
 \nonumber\\
 &\tL^\pm_{21} \,\tL^\pm_{22} = 0
 &
 &\tL^\pm_{22} \,\tL^\pm_{12} =0
\end{alignat}
the relations (\ref{eq:RL+L-}) become
\begin{alignat}{2}
 \label{eq:RtLtL2}
 &
 [\tL^+_{11}, \tL^-_{11}] = 0
 &\qquad\qquad
 &
 \tL^-_{21} \,\tL^+_{11} = \tL^+_{21} \,\tL^-_{11}
 \nonumber\\
 &
 \tL^-_{11} \,\tL^+_{12} = \tL^+_{11} \,\tL^-_{12}
 &\qquad\qquad
 &
 \tL^-_{21} \,\tL^+_{12} = \tL^+_{21} \,\tL^-_{12}
 \nonumber\\
 &
 \tL^-_{11} \,\tL^+_{21} = \tL^+_{21} \,\tL^-_{22}
 &\qquad\qquad
 &
 \tL^-_{21} \,\tL^+_{21} = - \tL^+_{11} \,\tL^-_{22}
 \nonumber\\
 &
 \tL^-_{11} \,\tL^+_{22} = \tL^+_{21} \,\tL^-_{21}
 &\qquad\qquad
 &
 \tL^-_{21} \,\tL^+_{22} = - \tL^+_{11} \,\tL^-_{21}
 \nonumber\\[3mm]
 &
 \tL^-_{12} \,\tL^+_{11} = - \tL^+_{22} \,\tL^-_{12}
 &\qquad\qquad
 &
 \tL^-_{22} \,\tL^+_{11} = \tL^+_{12} \,\tL^-_{12}
 \nonumber\\
 &
 \tL^-_{12} \,\tL^+_{12} = - \tL^+_{22} \,\tL^-_{11}
 &\qquad\qquad
 &
 \tL^-_{22} \,\tL^+_{12} = \tL^+_{12} \,\tL^-_{11}
 \nonumber\\
 &
 \tL^-_{12} \,\tL^+_{21} = \tL^+_{12} \,\tL^-_{21}
 &\qquad\qquad
 &
 \tL^-_{22} \,\tL^+_{21} = \tL^+_{22} \,\tL^-_{21}
 \nonumber\\
 &
 \tL^-_{12} \,\tL^+_{22} = \tL^+_{12} \,\tL^-_{22}
 &\qquad\qquad
 &
 [\tL^+_{22}, \tL^-_{22}] = 0
\end{alignat}

We would like to introduce a basis for the FRT dual algebra. We
need the following notation:
\begin{eqnarray} \label{newbasis}
&&  F_n(k_i;l_i) \equiv \prod^n_{i=1}\tilde{L}^{+
    k_i}_{11} \tilde{L}^+_{12} \tilde{L}^{+
    l_i}_{22}\tilde{L}^{+}_{21}\ , \qquad n\geq 1, \nn\\
&&  G_n(l_i; k_i) \equiv \prod^n_{i=1}
  \tilde{L}^{+ l_i}_{22} \tilde{L}^+_{21} \tilde{L}^{+ k_i}_{11}
  \tilde{L}^+_{12}\ , \qquad n\geq 1, \nn\\
&& F_0(k_i;l_i) \equiv 1; \ \ \ G_0(l_i;k_i) \equiv 1 \ .
\end{eqnarray}

The basis elements of the algebra generated by the $\tL^+$'s are:
\begin{eqnarray}
  \label{eq:basisB+}
  F_n(k_i;l_i)\tilde{L}^{+k_{n+1}}_{11}\ , \qquad  F_{n-1}(k_i;l_i)
  \tilde{L}^{+k_n}_{11}\tilde{L}^+_{12}\tilde{L}^{+l_n}_{22}\ ,
  \nonumber\\
  G_n(l_i;k_i) \tilde{L}^{+l_{n+1}}_{22}\ ,  \qquad
  G_{n-1}(l_i;k_i)\tilde{L}^{+l_n}_{22}\tilde{L}^+_{21}
  \tilde{L}^{+k_n}_{11}\ .
\end{eqnarray}
Defining also $  K_n = \sum^n_{i=1} k_i, \ \ \ L_n = \sum^n_{i=1}
l_i $ the actions of generators $\tL^-$ on the basis elements are,
e.g.,
\begin{eqnarray}
  \tilde{L}^-_{11} F_n(k_i;l_i)&=& F_{n-1}(k_1
  +1,k_i;l_i) \tilde{L}^{+ k_n}_{11} \tilde{L}^+_{12} \tilde{L}^{+
    l_n}_{22} \tilde{L}^-_{21}\ ,  \nonumber \\ \tilde{L}^-_{12}
  F_n(k_i;l_i) &=& (-1)^{K_n + L_n +1} G_{n-1}(k_1+1,k_i;l_i)
  \tilde{L}^{+ k_n}_{22} \tilde{L}^+_{21}
  \tilde{L}^{+ l_n}_{11} \tilde{l}^-_{22}\ , \nonumber\\
  \tilde{L}^-_{21} F_n(k_i;l_i) &=& G_n(0,..,l_{n-1};k_i) \tilde{L}^{+
    l_n}_{22} \tilde{L}^-_{21}\ ,  \nonumber \\
  \tilde{L}^-_{22}
  F_n(k_i,;l_i) &=& (-1)^{K_n + L_n}F_n(0,..,l_{n-1};k_i) \tilde{L}^{+
    l_n}_{11} \tilde{L}^-_{22}\ ,
\end{eqnarray}
\begin{eqnarray}
  \tilde{L}^-_{11} G_n(l_i;k_i) &=&
  (-1)^{K_n+L_n}G_n(0, ..,k_{n-1};l_i) \tilde{L}^{+
    k_n}_{22}\tilde{L}^-_{11}\ , \nonumber\\
  \tilde{L}^-_{12}
    G_n(l_i;k_i) &=&  F_n(0,..,k_{n-1};l_i) \tilde{L}^{+
    k_n}_{11}\tilde{l}^-_{12}\ , \nn\\
  \tilde{L}^-_{21} G_n(l_i;k_i) &=& (-1)^{K_n+L_n+1}
  F_{n-1}(l_1+1,l_i;k_i)
  \tilde{L}^{+ l_n}_{11} \tilde{L}^+_{12}
  \tilde{L}^{+k_n}_{22}\tilde{L}^-_{11}\ , \nonumber \\
  \tilde{L}^-_{22} G_n(l_i;k_i) &=&
  G_{n-1}(l_1+1,l_i;k_i)\tilde{L}^{+
    l_n}_{22}\tilde{L}^+_{21}\tilde{L}^{+k_n}_{11} \tilde{L}^-_{12}\ .
\end{eqnarray}
These equations allow one to order the $\tL^-$ with respect to the
$\tL^+$. For the $\tL^-$ among themselves, there exists a basis
similar to (\ref{eq:basisB+}). Thus, this basis gives the class of representations
when $s03_F$ acts on itself.

\subsection{Affine $s03_F$}
\label{sect:affine}

$L^\pm(x)$ are now matrices of operators $L^\pm_{ij}(x)$ ($i,j=1,2$)
satisfying the relations
\begin{eqnarray}
 \label{eq:RLLaffine}
 R^+(x_1/x_2) \; L^+_1(x_1) \, L^+_2(x_2)
 = L^+_2(x_2) \, L^+_1(x_1) \; R^+(x_1/x_2) \nn\\
 R^+(x_1/x_2) \; L^-_1(x_1) \, L^-_2(x_2)
 = L^-_2(x_2) \, L^-_1(x_1) \; R^+(x_1/x_2) \nn\\
 R^+(x_1/x_2) \; L^+_1(x_1) \, L^-_2(x_2)
 = L^-_2(x_2) \, L^+_1(x_1) \; R^+(x_1/x_2)
\end{eqnarray}
i.e.
\begin{eqnarray}
 \label{eq:RLLaffine2}
 &&(x_1+x_2)
 \left( L^+_{ab}(x_1) L^+_{cd}(x_2) - L^+_{ab}(x_2) L^+_{cd}(x_1)\right)
 +\nonumber\\
 &&\qquad
 + \theta_a (x_2-x_1) L^+_{\bar a b}(x_1) L^+_{\bar c d}(x_2)
 + \theta_b (x_2-x_1) L^+_{a \bar b}(x_1) L^+_{c \bar d}(x_2) = 0 \qquad
\nn \\[2mm]
 &&(x_1+x_2)
 \left( L^-_{ab}(x_1) L^-_{cd}(x_2) - L^-_{ab}(x_2) L^-_{cd}(x_1)\right)
 +\nonumber\\
 &&\qquad
 + \theta_a (x_2-x_1) L^-_{\bar a b}(x_1) L^-_{\bar c d}(x_2)
 + \theta_b (x_2-x_1) L^-_{a \bar b}(x_1) L^-_{c \bar d}(x_2) = 0 \qquad
\nn \\[2mm]
 &&(x_1+x_2)
 \left( L^+_{ab}(x_1) L^-_{cd}(x_2) - L^-_{ab}(x_2) L^+_{cd}(x_1)\right)
 +\nonumber\\
 &&\qquad
 + \theta_a (x_2-x_1) L^+_{\bar a b}(x_1) L^-_{\bar c d}(x_2)
 + \theta_b (x_2-x_1) L^-_{a \bar b}(x_1) L^+_{c \bar d}(x_2) = 0 \qquad
\end{eqnarray}

In particular,
\begin{eqnarray}
 \left[ L_{12}^+(x_1) , L_{12}^-(x_2) \right] -
 \left[ L_{21}^+(x_1) , L_{21}^-(x_2) \right] = 0 \nn\\
 \left[ L_{11}^+(x_1) , L_{11}^-(x_2) \right] +
 \left[ L_{22}^+(x_1) , L_{22}^-(x_2) \right] = 0
\end{eqnarray}

\section{Alternative parametrisation}
\setcounter{equation}{0}

Another parametrisation of the braid matrix is:

\begin{equation} \hat{R}(z) = {{(1+z)\hat{R} + (1-z)\hat{R}^{-1}}\over (2(1+z^2))^{1/2}}
= {1\over \sqrt{1+z^2}}\ \left(
    \begin{array}{cccc}
      1 & 0 & 0 & z \cr
      0 & 1 & -z & 0 \cr
      0 & z & 1 & 0 \cr
      -z & 0 & 0 & 1
    \end{array}
  \right)
\end{equation}

Some advantages of this parametrisation are:
\begin{equation}  \hat{R}(\pm 1) = \hat{R}^{\pm 1} \end{equation}
and
\begin{equation} \hat{R}_{12}(z'')\hat{R}_{23}(z)\hat{R}_{12}(z') = \hat{R}_{23}(z')\hat{R}_{12}(z)\hat{R}_{23}(z'')
\end{equation}
where
\begin{equation} z'' = {{z-z'}\over {1-zz'}} \end{equation}

A Baxterisation for $L^\pm$ is:
\begin{equation} L(z) = {{(1+z)L^{+} + (1-z)L^{-}}\over (2(1+z^2))^{1/2}} \end{equation}
here
\begin{equation} L(\pm 1) = L^{\pm} \end{equation}
and
\begin{equation}  \hat{R}(z'')L_2(z)L_1(z') = L_2(z')L_1(z)\hat{R}(z'') \end{equation}

If we accept the convention $z'' =1$ for $z = \pm 1$ and $z' = \pm
1$ we can reproduce the formulae for $L^{\pm}$
(\ref{eq:RLL},\ref{eq:pybel},\ref{eq:RLLij},\ref{eq:RL+L-}).

 For ${\hat R}$ matrices satisfying a minimal quadratic
equation (the first equation of \eqref{eq:hRprop}
being an example) there exist two
possibilities of defining co-products of L. Here they correspond to A)  and
B) below:\nl
A)
\begin{equation} \delta\, L_{ij}(z) = \Sigma_k L_{ik}(z) \otimes L_{kj}(z) \end{equation}
B)
\begin{equation} \tilde{\delta}\, L_{ij}(z) = {1\over (2(1+z^2)^{1/2}}((1-z)\delta L^+_{ij} + (1+z) \delta L^-_{ij})
\end{equation}
where
\begin{equation} L (z)= \left( \begin{array}{cc} L_{11}(z) & L_{12} (z)\cr L_{21}(z) & L_{22}(z) \end{array}
\right) \end{equation}

Both types of coproducts satisfy the equations
\begin{equation} \hat{R}(z'')(\delta L(z))_2 (\delta L(z'))_1 = (\delta L(z'))_2 (\delta L(z))_1 \hat{R}(z'') \end{equation}
and exactly the same for $\tilde{\delta}L$.

A two-dimensional representation for the algebra generated by the
$L$-operators is provided by the $R$-matrix itself, setting $\pi(L^+)
= R_{21}$, $\pi(L^-) = R^{-1}$ (see \cite{FRT,AC0401}).
Thus, if for $L$ is used the fundamental representation
\begin{equation}\label{lfund} L(z) = \left( \begin{array}{cc} L_{11} (z)& L_{12}(z) \cr L_{21}(z) & L_{22}(z) \end{array}
 \right) = \hat{R}(z)P =
{1\over  \sqrt{1+z^2}}\left( \begin{array}{cccc} 1 & 0 & 0 & z
\cr 0 & -z & 1 & 0 \cr 0 & 1 & z & 0 \cr -z & 0 & 0 & 1
\end{array} \right) \end{equation}
we have the following explicit forms
\begin{eqnarray}\label{ldelta}
& \delta L_{11}(z) = \left(  \begin{array}{cccc}
 1 & 0 & 0 & z \cr
 0 & -z & -z^2 & 0 \cr
 0 & 1 & -z & 0 \cr
 -z & 0 & 0 & z^2
 \end{array} \right), \ \ \ &\tilde{\delta} L_{11}(z) = \left(
 \begin{array}{cccc}
 1 & 0 & 0 & z \cr
 0 & -z & -1 & 0 \cr
 0 & 1 & -z & 0 \cr
 -z & 0 & 0 & 1
 \end{array} \right) \nn\\
&  \delta L_{12}(z) = \left(  \begin{array}{cccc}
 0 & z & z^2 & 0 \cr
 1 & 0 & 0 & z \cr
 z & 0 &  &  -z^2 \cr
 0 & 1 & -z & 0
 \end{array} \right), \ \ \ &\tilde{\delta} L_{12}(z) = \left(
 \begin{array}{cccc}
 0 & z & 1 & 0 \cr
 1 & 0 & 0 & z \cr
 z & 0 & 0 & -1 \cr
 0 & 1 & -z & 0
 \end{array} \right) \nn\\
&  \delta L_{21}(z) = \left(  \begin{array}{cccc}
 0 & z & 1 & 0 \cr
 -z^2 & 0 & 0 & -z \cr
 -z & 0 & 0 &  1 \cr
 0 & z^2 & -z & 0
 \end{array} \right), \ \ \ &\tilde{\delta} L_{21}(z) = \left(
 \begin{array}{cccc}
 0 & z & 1 & 0 \cr
 -1 & 0 & 0 & -z \cr
 -z & 0 & 0 & 1 \cr
 0 & 1 & -z & 0
 \end{array} \right) \nn\\
&  \delta L_{22}(z) = \left(  \begin{array}{cccc}
 z^2 & 0 & 0 & z \cr
 0 & z & 1 & 0 \cr
 0 & -z^2 & z &  0 \cr
 -z & 0 & 0 & 1
 \end{array} \right), \ \ \ &\tilde{\delta} L_{22}(z) = \left(
 \begin{array}{cccc}
 1 & 0 & 0 & z \cr
 0 & z & 1 & 0 \cr
 0 & -1 & z & 0 \cr
 -z & 0 & 0 & 1
 \end{array} \right)
\end{eqnarray}
 These two sets coincide for $z=\pm1$ but except for
these limits they  can be shown to be inequivalent.

We use the diagonaliser $M$ (cf. \cite{Chakra3}):
\begin{equation}
  M = \frac1{\sqrt2}
 \left(
 \begin{array}{cccc}
 1 & 0 & 0 & i \cr
 0 & 1 & -i & 0 \cr
 0 & -i & 1 & 0 \cr
 i & 0 & 0 & 1
 \end{array}
 \right) \ .
\end{equation}
to define the matrices $X$
\begin{equation} M L_2(z) L_1(z') M^{-1} = X(z,z') \ , \qquad
 M L_2(z') L_1(z) M^{-1} = X'(z',z) \end{equation}
so that
\begin{equation} X_{ij}(z,z') = X'_{ij}(z',z)    \ \ \ i,j = 1,2,3,4. \end{equation}

Denote $L(z) = L$ and $L(z') = L'$.
Then we have for the components of $X$:
\begin{eqnarray}\label{xcomp}
&& \left( \begin{array}{c} X_{11} \cr X_{44} \end{array} \right) = (L_{11}L'_{11} +
L_{22}L'_{22}) \pm i (L_{21} L'_{21} - L_{12}L'_{12}) \nn\\
&& \left( \begin{array}{c} X_{12} \cr X_{43} \end{array} \right) = \pm (L_{12}L'_{11} -
L_{21}L'_{22}) + i (L_{22} L'_{21} + L_{11}L'_{12}) \nn\\
&& \left( \begin{array}{c} X_{21} \cr X_{34} \end{array} \right) = \pm (L_{21}L'_{11} -
L_{12}L'_{22}) - i (L_{11} L'_{21} + L_{22}L'_{12}) \nn\\
&& \left( \begin{array}{c} X_{22} \cr X_{33} \end{array} \right) = (L_{22}L'_{11} +
L_{11}L'_{22}) \mp i (L_{12} L'_{21} - L_{21}L'_{12}) \nn\\
&& \left( \begin{array}{c} X_{13} \cr X_{42} \end{array} \right) = \pm (L_{11}L'_{12} -
L_{22}L'_{21}) + i (L_{12} L'_{11} + L_{21}L'_{22}) \nn\\
&& \left( \begin{array}{c} X_{23} \cr X_{32} \end{array} \right) =  (L_{12}L'_{21} +
L_{21}L'_{12}) \pm i (L_{22} L'_{11} - L_{11}L'_{22}) \nn\\
&& \left( \begin{array}{c} X_{14} \cr X_{41} \end{array} \right) =  (L_{21}L'_{21} +
L_{12}L'_{12}) \mp i (L_{11} L'_{11} - L_{22}L'_{22}) \nn\\
&& \left( \begin{array}{c} X_{24} \cr X_{31} \end{array} \right) =  \mp (L_{11}L'_{21} -
L_{22}L'_{12}) - i (L_{21} L'_{11} + L_{12}L'_{22}). \end{eqnarray}

Then the Yang-Baxter equation (3.4) reads:
\begin{equation} (\hat{R}(z''))_\alpha X = X' (\hat{R}(z''))_\alpha \end{equation}
which gives for the $\hat{R}L_2 L_1$ relations the following explicit formulae
\begin{eqnarray}
&& X_{ij}(z,z') = X'_{ij}(z',z)\ , \qquad  (ij)= (11,12,21,22;33,34,43,44)\ , \\
&& ((1-zz') -i(z-z')) X_{ij} = ((1-zz') +i(z-z')) X'_{ij}\ , \qquad   (ij) = (13,14,23,24)\ , \nn\\
&&((1-zz') +i(z-z')) X_{ij} = ((1-zz')-i(z-z')) X'_{ij}\ , \qquad  (ij) = (31,41,32,42)
  \nn
 \end{eqnarray}

\vspace{1cm}

Having in mind \eqref{xcomp}  the general structure of the $X_{ij}$ (up to normalisation factors) is:
\begin{eqnarray}\label{xxcomp}
&& X_{ij} = (1+zz')A_{ij} + (1+zz')B_{ij} + (z+z')C_{ij}\ ,
\ \ \  (ij) = (11,12,21,22;33,34,43,44)\nn\\
&& X_{ij} = ((1-zz' +i(z-z')) Q_{ij}\ , \ \ \  (ij) = (13,14,23,24)\nn\\
&& X_{ij} = ((1-zz' -i(z-z')) N_{ij}\ , \ \ \  (ij) = (31,41,32,42) \end{eqnarray}
where the matrices $\{A_{ij}, B_{ij}, C_{ij}, Q_{ij}, N_{ij} \}$
do not depend on $z$. However they can not be arbitrary, but should
be compatible with the definitions (\ref{xcomp}) for $X_{ij}$.

{\it An example:}\ \ \
If we take the fundamental representation for $L(z)$ (\ref{lfund}) we obtain
\begin{eqnarray}\label{xxxcomp} &&
\left( \begin{array}{cc} X_{11} & X_{12} \cr X_{21} & X_{22} \end{array} \right)
= ((1+zz') - i(z+z'))
\left( \begin{array}{cccc} 1 & 0 & 0 & -1 \cr 0 & 1 & 1 & 0 \cr 0
& -i & i & 0 \cr -i & 0 & 0 & -i \end{array} \right)\nn\\ &&
 \left( \begin{array}{cc} X_{33} & X_{34} \cr X_{43} & X_{44} \end{array} \right)
 = ((1+zz') + i(z+z'))
\left( \begin{array}{cccc} -i & 0 & 0 & -i \cr 0 & i & -i & 0 \cr
0 & 1 & 1 & 0 \cr -1 & 0 & 0 & 1 \end{array} \right) \nn\\ &&
\left( \begin{array}{cc} X_{13} & X_{14} \cr X_{23} & X_{24} \end{array} \right) = i((1-zz') + i(z-z'))
\left( \begin{array}{cccc} 0 & 1 & 1 & 0 \cr 1 & 0 & 0 & -1 \cr -i & 0 & 0 & i \cr 0 & -i & -i & 0 \end{array}
\right) \nn\\ &&
 \left( \begin{array}{cc} X_{31} & X_{32} \cr X_{41} & X_{42} \end{array} \right) = -i((1-zz') + i(z-z'))
\left( \begin{array}{cccc} 0 & i & i & 0 \cr -i & 0 & 0 & i \cr 1 & 0 & 0 & 1 \cr 0 & -1 & 1 & 0 \end{array} \right)
\end{eqnarray}
So in that case:
\begin{equation} B_{ij} = 0 , \ \ \ C_{ij} = \mp i A_{ij} \ .
\end{equation}

A class of representations for arbitrary dimensions can be considered.
Based on formulae \eqref{ldelta}  and subsequent use of the diagonaliser $M$ we have:
\begin{equation} M(\tilde{\delta} L_{ij})M^{-1} = \left( \begin{array}{cc} U_{ij} & 0 \cr 0 & D_{ij} \end{array} \right)
\end{equation}
where (neglecting normalisation factors):
\begin{equation}\label{ud}
\begin{array}{l}
 (U_{11}, U_{22}, U_{12}, U_{21}) = (1-iz) \left( \left( \begin{array}{cc} 1 & 0 \cr 0 & -i \end{array}
\right) , \left( \begin{array}{cc} 1 & 0 \cr 0 & i \end{array}
\right) , \left( \begin{array}{cc} 0 & i \cr 1 & 0 \end{array}
\right) , \left( \begin{array}{cc} 0 & i \cr -1 & 0 \end{array}
\right)  \right)\ , \cr \cr 
  (D_{11}, D_{22}, D_{12}, D_{21}) = i(1+iz) \left( \left( \begin{array}{cc} 1 & 0 \cr 0 & -i \end{array}
\right) , \left( \begin{array}{cc} -1 & 0 \cr 0 & -i \end{array}
\right) , \left( \begin{array}{cc} 0 & i \cr 1 & 0 \end{array}
\right) , \left( \begin{array}{cc} 0 & -i \cr 1 & 0 \end{array}
\right)  \right)\ .  \end{array}
\end{equation}
For $z = \pm 1$ \ \  $U^{\pm}_{ij}, D^{\pm}_{ij}$  give a
particular class of complex $2 \times 2$ representations.

In view of \eqref{ud}, where  there are  factors $(1 \pm iz)$ appearing, we try the Ansatz:
\begin{equation} L_{ij} = (1+kz) \hat{L}_{ij} \end{equation}
where $\hat{L}_{ij}$ is $z$-independent. Then each $X_{ij}$ would
be proportional to $(1+kz)(1+kz')$ and  relations
\eqref{xxcomp} are satisfied with:
\begin{equation}\label{ansa}  (A_{ij} - B_{ij}) = k^2(A_{ij} + B_{ij}) = k C_{ij} \ , \qquad
  Q_{ij} = 0, \quad N_{ij} = 0\ . \end{equation}

Now consider as an example the $ 3 \times 3$ case:
\begin{eqnarray}\label{thr} & L_{11} = (1+kz) \left( \begin{array}{ccc} a & 0 & 0 \cr  0 & b & 0 \cr 0 & 0 & c \end{array}
\right)\ , \ \ \ &L_{22} = \pm (1+kz) \left( \begin{array}{ccc} a & 0 & 0 \cr
 0 & -b & 0 \cr 0 & 0 & c \end{array} \right)\ , \nn\\
& L_{12} = (1+kz) \left( \begin{array}{ccc} 0 & d & 0 \cr  e & 0 & f \cr 0 & g & 0 \end{array} \right)
\ , \ \ \
&L_{21} = \pm (1+kz) \left( \begin{array}{ccc} 0 & d & 0 \cr  -e & 0 & -f \cr 0 & g & 0 \end{array} \right)
\end{eqnarray}
It can be seen that equations \eqref{ansa} are satisfied.

The analogue of \eqref{thr} for the general $n \times n$ case is:
\begin{eqnarray}\label{ntn} && (L_{11})_{mn} = (1+kz)\,\delta_{m,n}\, a_m\ ,
\nn\\  &&  (L_{22})_{mn} = \epsilon\,(1+kz)\, \delta_{m,n} (-1)^{m-1}\, a_m\ ,
\nn\\  &&    (L_{12})_{mn} = (1+kz)\,(\delta_{n,m+1} u_m + \delta_{n,m-1} v_m)\ ,
\nn\\  &&  (L_{21})_{mn} = \epsilon\,(1+kz)\, (\delta_{n,m+1} (-1)^{m+1}u_m + \delta_{n,m-1} (-1)^{m-1}v_m)
\end{eqnarray}


\section{Finite dimensional representations}
\label{sect:RRR}
\setcounter{equation}{0}
\subsection{Representations on S03}

Here we shall study the representations of ~$s03_F$~ obtained by the
use of its right regular action (RRA) on the dual bialgebra ~$S03$.
 The RRA is defined as follows:
\begin{equation}
 \label{eq:RRRdef}
 \pi_R(L^\pm_{ij}) f = f^{(1)} <L^\pm_{ij}, f_{(2)}>
\end{equation}
where we use Sweedler's notation for the co-product: ~$\d(f) ~=~
f_{(1)} \otimes f_{(2)}\,$.   More explicitly, for the generators of
~$s03_F$~ we have:
\begin{alignat}{2}
 \label{eq:RRRabcd}
 &
 \pi_R(L^\pm_{11})
 \left(
 \begin{array}{cc}
 a & b \cr
 c & d
 \end{array}
 \right)
 =
 \left(
 \begin{array}{cc}
 a & \mp b \cr
 c & \mp d
 \end{array}
 \right)
 &\qquad
 &
 \pi_R(L^\pm_{12})
 \left(
 \begin{array}{cc}
 a & b \cr
 c & d
 \end{array}
 \right)
 =
 \left(
 \begin{array}{cc}
 b & \pm a \cr
 d & \pm c
 \end{array}
 \right)
 \nonumber\\
 &
 \pi_R(L^\pm_{21})
 \left(
 \begin{array}{cc}
 a & b \cr
 c & d
 \end{array}
 \right)
 =
 \left(
 \begin{array}{cc}
 \mp b & a \cr
 \mp d & c
 \end{array}
 \right)
 &\qquad
 &
 \pi_R(L^\pm_{22})
 \left(
 \begin{array}{cc}
 a & b \cr
 c & d
 \end{array}
 \right)
 =
 \left(
 \begin{array}{cc}
 \pm a & b \cr
 \pm c & d
 \end{array}
 \right)
\end{alignat}
Obviously, the above representation is the direct sum of two
equivalent 
two-dimensional irreps (with vector spaces spanned by $\{a,b\}$
and $\{c,d\}$, respectively) such that the representation matrices 
(acting on $(a,b)$ or $(c,d)$ from the right) are given by:
\begin{alignat}{2}
 \label{eq:RRRlinear}
 &
 \pi_R(L^\pm_{11})
 =
 \left(
 \begin{array}{cc}
 1 & 0 \cr
 0 & \mp 1
 \end{array}
 \right)
 &\qquad\qquad
 &
 \pi_R(L^\pm_{12})
 =
 \left(
 \begin{array}{cc}
 0 & \pm 1 \cr
 1 & 0
 \end{array}
 \right)
 \nonumber\\
 &
 \pi_R(L^\pm_{21})
 =
 \left(
 \begin{array}{cc}
 0 & 1 \cr
 \mp 1 & 0
 \end{array}
 \right)
 &\qquad
 &
 \pi_R(L^\pm_{22})
 =
 \left(
 \begin{array}{cc}
 \pm 1 & 0 \cr
 0 & 1
 \end{array}
 \right)
\end{alignat}

Further, we would like to consider polynomials in the elements $a,b,c,d$
of degree $N>1$. Superficially, for fixed $N$ such polynomials would
span a vector space of dimension $4^N$, however, 
due to the relations \eqref{eq:S03rel} such polynomials actually 
span a vector space of dimension $2^{N+1}$. Explicitly, these vector spaces 
are spanned by:
\eqn{monobas}
\{a,b\}^{\otimes^N} \ , 
\quad  \{c,d\}\otimes \{a,b\}^{\otimes^{N-1}} \ , 
\end{equation}
These vector spaces split into irreducible 
representations for which the most suitable bases are complex linear
combinations of the above. We give explicitly some examples of small $N$
and then formulate a general statement.

For $N=2$ from the vector space of dimension 8 one can
extract four two-dimensional irreducible representations (two
by two equivalent).
The representation spaces are spanned over the
elements $V^\epsilon_1 = a^2 + i\epsilon b^2, \ \ V^\epsilon_2 =
ab - i\epsilon ba$, where $\epsilon = \pm 1$ labels the two 
non-equivalent irreducible representations.  
Another two irreducible representations
are spanned over the elements $\tilde{V}^\epsilon_1 = ca
+i\epsilon db, \ \ \tilde{V}^\epsilon_2 = cb -i\epsilon da$, 
however, they are  equivalent to the first two for coinciding 
values of $\eps$. In a
matrix form the representations are as follows:
\begin{eqnarray}
&&\pi_R(L^\pm_{11}) = (1 \pm i\epsilon)\left(\begin{array}{cc} 1 & 0 \\
0 &  - i\epsilon \end{array} \right), \ \ \ \pi_R(L^\pm_{12}) = (1
\pm i\epsilon)\left(\begin{array}{cc}0 &  -i\epsilon \\
1 & 0 \end{array} \right), \nn\\
&&\pi_R(L^\pm_{21}) = (1 \pm i\epsilon)
\left(\begin{array}{cc} 0 & -i\epsilon \\
-1 & 0 \end{array} \right), \ \ \ \pi_R(L^\pm_{22}) =
(1 \pm i\epsilon)
\left(\begin{array}{cc}1  & 0 \\
0 & i\epsilon \end{array} \right).
\end{eqnarray}
(In these formulae the $\epsilon$ just indicates the the
representation concerned.)

Having these elements one can proceed further to construct all 
representations for any fixed $N$. 
For $N=3$ the overall vector space is 16-dimensional. 
We first consider the 8-dimensional vector space 
~$\{a,b\}^{\otimes^3}$~ for which the convenient basis is:
\begin{equation}
U_1^\epsilon = aV_1^\epsilon, \ \ \ 
U_2^\epsilon = bV_2^\epsilon, \ \ \ 
U_3^\epsilon = aV_2^\epsilon, \ \ \ 
U_4^\epsilon = bV_1^\epsilon \ .
\end{equation}
These elements form two four-dimensional  irreducible
representations, labelled again by the index $\epsilon$. 
In  matrix form these representations can be written as:
\begin{eqnarray}\label{cubic} 
\pi_R(L^{\pm}_{11}) = (1\pm i\epsilon)\left( 
\begin{array}{cccc}
1  & \mp i\epsilon & 0 & 0 \\ -1 & \mp i\epsilon
 & 0 & 0 \\ 0 & 0 &  i\epsilon & \mp 1 \\
0 & 0 & -i \epsilon &  \mp 1 \end{array} \right), 
~\pi_R(L^{\pm}_{12}) = (1 \pm i\epsilon)\left(
\begin{array}{cccc} 0 & 0 &  -i\epsilon & \pm 1 \\ 0 & 0 & -i\epsilon
& \mp1 \\ 1 & \mp i\epsilon &
0 & 0 \\ 1  & \pm i\epsilon & 0 & 0 \end{array} \right) \\ \nonumber \\
\pi_R(L^{\pm}_{21}) = (1 \pm i\epsilon)\left(
\begin{array}{cccc}0 & 0 &  \mp i \epsilon & 1  \\ 0 & 0 &  \mp i\epsilon
& -1  \\ \mp 1  & i\epsilon & 0 & 0
\\ \mp 1  & - i\epsilon) & 0 & 0
\end{array} \right), 
~\pi_R(L^{\pm}_{22}) = (1\pm i\epsilon)\left( \begin{array}{cccc}
\pm 1 & - i\epsilon & 0 & 0 \\  \mp 1 & -i\epsilon
 & 0 & 0 \\ 0 & 0 &  \mp i\epsilon &  1 \\
0 & 0 & \pm i \epsilon & 1 \end{array} \right) \nn 
\end{eqnarray}
Clearly, these four-dimensional representations are irreducible. 
For the remaining 8-dimensional vector space coming from 
~$\{c,d\}\otimes \{a,b\}^{\otimes^2}$~  the convenient basis is:
\begin{equation} 
\tilde{U}_1^{\epsilon} = cV_1^\epsilon, \ \ \ 
\tilde{U}_2^{\epsilon} = dV_2^\epsilon, \ \ \ 
\tilde{U}_3^{\epsilon} = cV_2^\epsilon, \ \ \ 
\tilde{U}_4^{\epsilon} = dV_1^\epsilon \ .
\end{equation}
The transformation rules for the elements $\tilde{U^\eps}$ 
are the same as those for ${U^\eps}$ given in \eqref{cubic} 
for the same values of $\eps$. 
Thus, again we have four irreducible representations,
which are two by two equivalent.

For $N=4$ the overall vector space is 32-dimensional. 
Using the elements $V^\epsilon_i$ and
$\tilde{V}^\epsilon_i$ it can be split into the following 8
four-dimensional representations:
\begin{equation} \omega^{\epsilon,\epsilon}_{ij} = V^\epsilon_i
V^\epsilon_j, \ \ \ \hat{\omega}^{\epsilon,-\epsilon}_{ij} =
V^\epsilon_i V^{-\epsilon}_j , \ \ \ 
\tilde{\omega}^{\epsilon,\epsilon}_{ij} = \tilde{V}^\epsilon_i
V^\epsilon_j, \ \ \ \hat{\tilde{\omega}}^{\epsilon,-\epsilon}_{ij}
= \tilde{V}^\epsilon_i V^{-\epsilon}_j\ , ~~~\eps=\pm 1, 
\end{equation}
(four sets doubled by $\eps$), 
where the indices $ij$ enumerate the four elements of a representation.
In a matrix form the representation formulae for $\omega$ are:
\begin{eqnarray}\label{quart} 
\pi_R(L^\pm_{11}) = \pm 2 \left(\begin{array}{cccc} i\epsilon & 0
& 0 & -i\epsilon \\ 0 & -1 & -1 & 0 \\  0 & 1 & -1 & 0 \\
-i\epsilon &  0 &  0 & -i\epsilon \end{array} \right),
~\pi_R(L^\pm_{12}) = \pm 2 \left(\begin{array}{cccc} 0 & 1 & 1 & 0 \\
i\epsilon & 0 & 0 & -i\epsilon \\  i\epsilon & 0 & 0 & i\epsilon \\
0 & 1 & -1 & 0 \end{array} \right) \\ \nonumber \\
\pi_R(L^\pm_{21}) = \pm 2 \left(\begin{array}{cccc} 0 & 1
& 1 & 0 \\ -i\epsilon & 0 & 0 & i\epsilon \\ -i\epsilon & 0 & 0 & 
-i\epsilon \\
0 & 1 & -1 & 0 \end{array} \right),  
~\pi_R(L^\pm_{22}) = \pm
2 \left(\begin{array}{cccc} i\epsilon & 0
& 0 & -i\epsilon \\ 0 & 1 & 1 & 0 \\  0 & -1 & -1 &  \\
-i\epsilon & 0 & 0 & -i\epsilon \end{array} \right)\ .\nn 
\end{eqnarray}
To write this matrix formula we used the conventional ordering of
the elements -- $\omega_{11}, \omega_{12}, \omega_{21},
\omega_{22}$. The analogous formulae for the $\hat{\omega}$ read:
\begin{eqnarray}\label{quartb}
\pi_R(L^\pm_{11}) =  2 \left(\begin{array}{cccc} 1 & 0
& 0 & 1 \\ 0 & -i\epsilon & i\epsilon & 0 \\  0 & i\epsilon & i\epsilon & 0 \\
-1 &  0 &  0 & 1 \end{array} \right),
~\pi_R(L^\pm_{12}) =  2 \left(\begin{array}{cccc} 0 & i\epsilon &
 -i\epsilon & 0 \\
1 & 0 & 0 & 1 \\ 1 & 0 & 0 & -1 \\
0 & i\epsilon & i\epsilon & 0 \end{array} \right) \\ \nonumber \\
\pi_R(L^\pm_{21}) =  2 \left(\begin{array}{cccc} 0 & i\epsilon
& i\epsilon & 0 \\ -1 & 0 & 0 & -1 \\ -1 & 0 & 0 & 1 \\
0 & i\epsilon & i\epsilon & 0 \end{array} \right), 
~\pi_R(L^\pm_{22}) = \pm 2 \left(\begin{array}{cccc} 1 & 0
& 0 & 1 \\ 0 & i\epsilon 1 & -i\epsilon & 0 \\  0 & -i\epsilon 
& -i\epsilon & 0 \\
-1 & 0 & 0 & 1 \end{array} \right)\ .\nn
\end{eqnarray}
One can check that the four-dimensional 
representations given by \eqref{quart} or \eqref{quartb}
 are irreducible. The representations 
$\tilde{\omega}^{\eps,\eps}$, ($\tilde{\hat{\omega}}^{\eps,-\eps}$),   
are equivalent to ${\omega}^{\eps,\eps}$, 
($\hat{\omega}^{\eps,-\eps}$),  
for the respective value of $\eps$.

We have carried out explicitly also the cases $N=5,6$ and all 
these results lead us to the following general statements.

In general the bases of degree $N=2n$ can be written in the
form:
\begin{equation} \Omega^{\epsilon_1,\epsilon_2, ...
,\epsilon_n}_{i_1,i_2, ...,i_n} =
V^{\epsilon_1}_{i_1}V^{\epsilon_2}_{i_2}...V^{\epsilon_n}_{i_n}, \
\ \ \tilde{\Omega}^{\epsilon_1,\epsilon_2, ...
,\epsilon_n}_{i_1,i_2, ...,i_n} =
\tilde{V}^{\epsilon_1}_{i_1}V^{\epsilon_2}_{i_2}...V^{\epsilon_n}_{i_n}
\end{equation}
where the set of indices $\{\epsilon_1,\epsilon_2, ... ,\epsilon_n\}$
labels the $2^n$ representations, while the indices 
$\{i_1,i_2,...,i_n\}$ enumerate the $2^n$ elements within a 
given representation.

The bases of odd order $N=2n+1$ are
constructed multiplying the above even elements $\Omega$ 
from the left by $a$ and $b$, then by $c$ and $d$, 
the second batch of representations being equivalent 
to the first. 

Thus, we can formulate the following general statement:
\begin{proposition}
 \begin{itemize}
 \item Tensor products of $2^{2n}$ 2-dimensional representations of the type
 described in
 (\ref{eq:RRRlinear}) (constructed using the coproduct structure)
 decompose into sums of  $2^{n+1}$ $2^n$-dimensional 
irreducible representations. 
These are 2-by-2 equivalent, i.e., 
the number of non-equivalent irreducible representations is $2^n$. 
 \item
 Tensor products of $2^{2n+1}$ 2-dimensional representations of the type
 described in 
 (\ref{eq:RRRlinear}) decompose into sums of $2^{n+1}$ 
$2^{n+1}$-dimensional irreducible representations. 
These are 2-by-2 equivalent, i.e., 
the number of non-equivalent irreducible representations is $2^n$. 
 \end{itemize}
\end{proposition}


\subsection{Finite dimensional irreducible representations 
(other constructions)}
\label{sect:irreps}

\subsubsection{Generalities}
$\tL^+_{11}$ and $\tL^-_{11}$ commute, so they have a common
eigenvector $v_0$.
\begin{equation}
 \label{eq:v0}
 \tL^+_{11}v_0 = \lambda^+ v_0 \qquad \tL^-_{11}v_0 = \lambda^- v_0
\end{equation}

\nt\textbf{A.  } Let us first suppose $\lambda^+ \neq 0$ and
$\lambda^- \neq 0$. Then
\begin{eqnarray}
 \tL^+_{22} v_0 = \tL^+_{12} v_0 = 0
\nn \\
 \tL^-_{22} v_0 = \tL^-_{12} v_0 = 0
\end{eqnarray}

On the whole representation,
\begin{equation}
 \label{eq:ratio}
 \tL^-_{ij} = \frac{\lambda^-}{\lambda^+} \tL^+_{ij}
\end{equation}
Indeed,
\begin{eqnarray}
 \tL^-_{i1} v_0 &=& \frac{1}{\lambda^+} \tL^-_{i1} \tL^+_{11} v_0
 \ =\ \frac{1}{\lambda^+} \tL^+_{i1} \tL^-_{11} v_0
 \ =\ \frac{\lambda^-}{\lambda^+} \tL^+_{i1} v_0
\end{eqnarray}
and, by recursion,
\begin{eqnarray}
 \tL^-_{ij_1} \big(\tL^+_{j_1 j_2} \tL^+_{j_2 j_3}\cdots \tL^+_{j_n
 j_{n+1}} \big) \;v_0 &=&
 \tL^+_{ij_1} \tL^-_{j_1 j_2} \big(\tL^+_{j_2 j_3}\cdots \tL^+_{j_n
 j_{n+1}} \big) \;v_0 \nn\\
 &=&
 \frac{\lambda^-}{\lambda^+} \
 \tL^+_{ij_1} \;\tL^+_{j_1 j_2} \;\big(\tL^+_{j_2 j_3}\cdots \tL^+_{j_n
 j_{n+1}} \big) \;v_0
\end{eqnarray}
when consecutive indices coincide, whereas in the other case \be
 \tL^-_{ij_1} \big(\tL^+_{\bar{j_1} j_2} \tL^+_{j_2 j_3}\cdots \tL^+_{j_n
 j_{n+1}} \big) \;v_0 =
 \pm \tL^+_{\bar ij_1} \tL^-_{\bar{j_1} \bar{j_2}}
 \big(\tL^+_{j_2 j_3}\cdots \tL^+_{j_n
 j_{n+1}} \big) \;v_0 = 0
\end{equation}

\textbf{B. } Let us consider the case $\lambda^+ \neq 0$ and
$\lambda^- = 0$.  This case is not yet completely understood. Let us
mention that the regular representation on linear terms in $a$, $b$,
$c$ and $d$ described by (\ref{eq:RRRlinear}) falls in this case.

A particular class of such representations corresponds to the choice
\begin{equation}
 \label{eq:subcase}
 L^-_{ii} = L^+_{\bar i \bar i} \qquad L^-_{i\bar i} = xL^+_{\bar i i}
\end{equation}
i.e.
\begin{equation}
 \label{eq:subcaset}
 \tL^-_{11} = \tL^+_{22}\,, \qquad
 \tL^-_{12} = x\tL^+_{21}\,, \qquad
 \tL^-_{21} = - x\tL^+_{12}\,, \qquad
 \tL^-_{22} = - \tL^+_{11}
\end{equation}
leading to supplementary relations for $\tL^+$ (with respect to
(\ref{eq:RtLtL1}))
\begin{alignat}{2}
 &
 \tL^+_{12} \tL^+_{21} = - x^{-1} (\tL^+_{11})^2
 &\qquad&
 \tL^+_{21} \tL^+_{12} = - x^{-1} (\tL^+_{22})^2
 \nn\\
 &
 \tL^+_{11} \tL^+_{12} = - \tL^+_{12} \tL^+_{22}
 &\qquad&
 \tL^+_{21} \tL^+_{11} = - \tL^+_{22} \tL^+_{21}
\end{alignat}

\subsubsection{2-dim irreps}
Two dimensional representations fall again into three cases.
\begin{itemize}
\item Those on which $\tL^\pm v_0 = \lambda^\pm v_0$ with both
 $\lambda^+$, $\lambda^-$ non-zero. They are described by
\begin{alignat}{2}
 \label{eq:2d1}
 &
 \pi(\tL^+_{11})
 =
 \left(
 \begin{array}{cc}
 \lambda^+ & 0 \cr
 0 & 0
 \end{array}
 \right)
 &\qquad\qquad
 &
 \pi(\tL^+_{12})
 =
 \left(
 \begin{array}{cc}
 0 & x \cr
 0 & 0
 \end{array}
 \right)
 \nonumber\\
 &
 \pi(\tL^+_{21})
 =
 \left(
 \begin{array}{cc}
 0 & 0 \cr
 x & 0
 \end{array}
 \right)
 &\qquad
 &
 \pi(\tL^+_{22})
 =
 \left(
 \begin{array}{cc}
 \ 0 & 0 \cr
 0 & \mu^+
 \end{array}
 \right)
\end{alignat}
and
$ \pi(\tL^-_{ij}) = \frac{\lambda^-}{\lambda^+} \pi(\tL^+_{ij}) $.

\item
Those on which there exists $v_0$ such that $\tL^\pm v_0 = \lambda^\pm
v_0$ with $\lambda^+\neq 0$, $\lambda^-=0$.
\begin{alignat}{2}
 \label{eq:2d2}
 &
 \pi(\tL^+_{11})
 =
 \left(
 \begin{array}{cc}
 \lambda^+ & 0 \cr
 0 & 0
 \end{array}
 \right)
 &\qquad\qquad
 &
 \pi(\tL^+_{12})
 =
 \left(
 \begin{array}{cc}
 0 & x \lambda^+ \cr
 0 & 0
 \end{array}
 \right)
 \nonumber\\
 &
 \pi(\tL^+_{21})
 =
 \left(
 \begin{array}{cc}
 0 & 0 \cr
 x \lambda^+ & 0
 \end{array}
 \right)
 &\qquad
 &
 \pi(\tL^+_{22})
 =
 \left(
 \begin{array}{cc}
 \ 0 & 0 \cr
 0 & -\lambda^+
 \end{array}
 \right)
\end{alignat}
\begin{alignat}{2}
 \label{eq:2d3}
 &
 \pi(\tL^-_{11})
 =
 \left(
 \begin{array}{cc}
 0 & 0 \cr
 0 & \mu
 \end{array}
 \right)
 &\qquad\qquad
 &
 \pi(\tL^-_{12})
 =
 \left(
 \begin{array}{cc}
 0 & 0 \cr
 x^{-1} \mu & 0
 \end{array}
 \right)
 \nonumber\\
 &
 \pi(\tL^-_{21})
 =
 \left(
 \begin{array}{cc}
 0 & - x^{-1} \mu \cr
 0 & 0
 \end{array}
 \right)
 &\qquad
 &
 \pi(\tL^-_{22})
 =
 \left(
 \begin{array}{cc}
 \ \mu & 0 \cr
 0 & 0
 \end{array}
 \right)
\end{alignat}

\item
Those representations with only 0 eigenvalues for $\tL^\pm_{11}$.
In the case of Jordan form for $\tL^\pm_{11}$, the corresponding
representation can be proved not to be irreducible. Hence
\begin{equation}
 \label{eq:2d3a}
 \pi(\tL^+_{11})
 = 0 \;, \qquad\qquad
 \pi(\tL^+_{22})
 = 0
\end{equation}
Then $\tL^\pm_{12}$ and $\tL^\pm_{21}$ should be of the form
\begin{alignat}{2}
 \label{eq:2d4}
 &
 \pi(\tL^\pm_{12})
 =
 \left(
 \begin{array}{cc}
 0 & \ell_{12}^\pm \cr
 0 & 0
 \end{array}
 \right)
 &\qquad\qquad
 &
 \pi(\tL^\pm_{21})
 = a^\pm
 \left(
 \begin{array}{cc}
 1 & b \cr
 -b^{-1} & -1
 \end{array}
 \right)
\end{alignat}
with $a^+ \ell_{12}^- = a^- \ell_{12}^+$.

\end{itemize}

\subsubsection{Other irreps}
\label{sect:otherirrep}
\begin{example}
 Here is an example (in the $\tL$ basis) of a finite dimensional irreducible
 representation of arbitrary dimension $N_1+N_2$, where  $N_1$ and $N_2$
are two non-negative integers:
 \begin{eqnarray}
 \label{eq:reprexample}
 &&
 \pi(\tL_{11}) =
 \diag(\rho_1,\cdots,\rho_{N_1},\underbrace{0,\cdots,0}_{N_2})
 \qquad\qquad \rho_i\neq\rho_j \qmbox{\rm for} i\neq j \nn\\
 &&\pi(\tL_{22}) =
 \diag(\underbrace{0,\cdots,0}_{N_1},\lambda_1,\cdots,\lambda_{N_2})
 \qquad\qquad \lambda_i\neq\lambda_j \qmbox{\rm for} i\neq j
\nn \\
 &&\left( \pi(\tL_{12}) \right)_{ij} \neq 0 \qmbox{\rm iff}
 i\in\{1,\cdots,N_1\}, \quad
 j\in\{N_1+1,\cdots,N_1+N_2\} \qquad
\nn \\
 &&\left( \pi(\tL_{21}) \right)_{ij} \neq 0 \qmbox{\rm iff}
 i\in\{N_1+1,\cdots,N_1+N_2\}, \quad
 j\in\{1,\cdots,N_1\}\ .
 \end{eqnarray}
\end{example}

\section{Diagonalisation of $\hR$, fusion and evaluation representations}
\label{sect:fusion}
\setcounter{equation}{0}

Due to (\ref{eq:hRx}), if $\hR$ is diagonalisable with the matrix $M$,
then $\hR(x)$ will be diagonalisable with the same matrix $M$
independent of $x$. Actually
\begin{equation}
 \label{eq:M}
 M = \frac1{\sqrt2}
 \left(
 \begin{array}{cccc}
 1 & 0 & 0 & i \cr
 0 & 1 & -i & 0 \cr
 0 & -i & 1 & 0 \cr
 i & 0 & 0 & 1
 \end{array}
 \right) \ .
\end{equation}
is such that
$M \hR(x) M^{-1} = \frac{\displaystyle 1+i}{\displaystyle 2\sqrt2}
\ \diag(x-i,x-i,1-ix,1-ix)$.

Let $\mu_1(x)$, $\mu_2(x)$ be the eigenvalues of $\hR(x)$, then
\begin{equation}
 \label{eq:Pi12}
 \Pi^{(1)} \equiv \frac{\hR(x) -\mu_2(x)}{\mu_1(x) -\mu_2(x)}\ ,
 \qquad\qquad
 \Pi^{(2)} \equiv \frac{\hR(x) -\mu_1(x)}{\mu_2(x) -\mu_1(x)}
\end{equation}
are projectors (${\Pi^{(i)}}^2 = \Pi^{(i)}$ and $\Pi^{(1)}+\Pi^{(2)} = 1$)
on the eigenspaces of $\hR(x)$. They are independent of $x$.

Taking the representations $\pi(L)$ given by (\ref{eq:RRRlinear}),
\begin{eqnarray}
 \label{eq:piL}
 &&\left(
 \pi L^+_{ij}
 \right)_{m}^{p} = (R_{21})_{im}^{jp} = R_{mi}^{pj} \qquad \qquad
 \left(
 \pi L^-_{ij}
 \right)_{m}^{p} = (R^{-1})_{im}^{jp}
\end{eqnarray}

\begin{eqnarray}
 \label{eq:DeltaL+}
 (\pi\otimes\pi) \left(
 \delta(L^+_{ij})
 \right)_{mn}^{pq}
 &=& \pi \left( L^+_{ik} \right)_{m}^{p} \pi \left( L^+_{kj} \right)_{n}^{q}
 \nonumber\\
 &=& R_{mi}^{pk} \ R_{nk}^{qj}
 \nonumber\\
 &=& R_{21} \ R_{31} \qmbox{(formally)}
\end{eqnarray}
Using the Yang--Baxter equation
$\hR_{23} \ R_{21} \ R_{31} = R_{21} \ R_{31} \ \hR_{23}$, one has
\begin{equation}
 \label{eq:RDeltaL+}
 \hR_{23} \ (\pi\otimes\pi) \delta(L^+)
 =
 (\pi\otimes\pi) \Delta(L^+) \hR_{23}
\end{equation}

Similarly
\begin{eqnarray}
 \label{eq:DeltaL-}
 (\pi\otimes\pi) \left(
 \delta(L^-_{ij})
 \right)_{mn}^{pq}
 &=& \pi \left( L^-_{ik} \right)_{m}^{p} \pi \left( L^-_{kj} \right)_{n}^{q}
 \nonumber\\
 &=& (R^{-1})_{im}^{kp} \ (R^{-1})_{kn}^{jq}
 \nonumber\\
 &=& (R^{-1})_{21} \ (R^{-1})_{31} \qmbox{(formally)}
\end{eqnarray}
Using the Yang--Baxter equation
$\hR_{23} \ (R^{-1})_{12} \ (R^{-1})_{13} =
(R^{-1})_{12} \ (R^{-1})_{13} \ \hR_{23}$, one has also
\begin{equation}
 \label{eq:RDeltaL-}
 \hR_{23} \ (\pi\otimes\pi) \delta(L^-)
 =
 (\pi\otimes\pi) \delta(L^-) \hR_{23}
\end{equation}

Hence
\begin{equation}
 \label{eq:RDeltaL2}
 \left[ \Pi^{(i)} , (\pi\otimes\pi) \delta(L^\pm) \right] = 0
\end{equation}
so that the eigenstates of $\hR$ are left invariant by the tensor
product of the fundamental representation.

We turn now to evaluation representations.
Noting that the characteristic polynomial of $\hR$ is of degree
two, we can define an evaluation representation by
\begin{equation}
 \label{eq:eval}
 L^+(x) = x^{-1} L^+ + L^- \ , \qquad
 L^-(x) = L^+ + x L^- \ .
\end{equation}
If $L^\pm$ are representations of the relations (\ref{eq:RLL})
then $L^\pm(x)$ are representations of (\ref{eq:RLLaffine2}).

Using the fact that $L^+$ and $L^-$ become identical on the tensor
product of the 2-dimensional fundamental representation given by
(\ref{eq:RRRlinear}), it is straightforward to see that the $x$
dependence completely factorises out for the corresponding
evaluation representation.

\section{Possible applications}
\label{sect:phys}
\setcounter{equation}{0}

We repeat the Baxterised $R$-matrix of S03 \eqref{eq:R_S03x} in order
to introduce necessary notation:
\begin{equation}
 \label{eq:baxtS03}
 R(u) =
 \left(
 \begin{array}{cccc}
 a_1(u) & 0 & 0 & d_1(u) \cr
 0 & b_1(u) & c_1(u) & 0 \cr
 0 & c_2(u) & b_2(u) & 0 \cr
 d_2(u) & 0 & 0 & a_2(u)
 \end{array}
 \right) = \frac{1}{\sqrt{2u}}
 \left(
 \begin{array}{cccc}
 u+1 & 0 & 0 & 1-u \cr
 0 & u+1 & u-1 & 0 \cr
 0 & 1-u & u+1 & 0 \cr
 u-1 & 0 & 0 & u+1
 \end{array}
 \right) \ .
\end{equation}
It satisfies the Yang--Baxter equation with spectral parameter
\begin{equation}
 \label{eq:YBEspec}
 R_{12}(z)R_{13}(zw)R_{23}(w)=R_{23}(w)R_{13}(zw)R_{12}(z)
\end{equation}

In this section we use some ingredients of the quantum inverse scattering method
\cite{Skl1} (for a book exposition, see
\cite{Kor}), however, we would not be able to follow it
throughout, due to the peculiarities of our exotic algebra, and,
on the other hand, we are able to use some simple procedures,
which work just in our situation.

\subsection{An exotic eight-vertex model}
\label{sect:8vertex}

An integrable vertex model can be constructed in the following way:

The entries of $R$ are interpreted as the Boltzmann weight of a
statistical model. The Yang--Baxter equation for $R$ leads to a kind
of star-triangle equation (in Baxter's terms) for the weights of the
model.

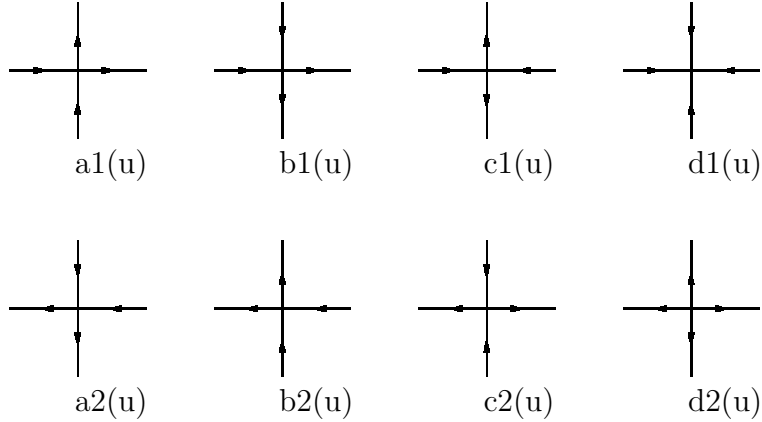
\begin{figure}[htbp]
 \centering
\setlength{\unitlength}{0.001 cm}

\begingroup\makeatletter\ifx\SetFigFont\undefined
\gdef\SetFigFont#1#2#3#4#5{
 \reset@font\fontsize{#1}{#2pt}
 \fontfamily{#3}\fontseries{#4}\fontshape{#5}
 \selectfont}
\fi\endgroup
{\renewcommand{\dashlinestretch}{30}
\begin{picture}(9924,5481)(0,-10)
\drawline(6312,4554)(6312,4104)
\drawline(6282.000,4224.000)(6312.000,4104.000)(6342.000,4224.000)
\drawline(6312,4104)(6312,3654)
\drawline(5412,4554)(5862,4554)
\drawline(5742.000,4524.000)(5862.000,4554.000)(5742.000,4584.000)
\drawline(5862,4554)(6312,4554)
\drawline(9042.000,3984.000)(9012.000,4104.000)(8982.000,3984.000)
\drawline(9012,4104)(9012,3654)
\drawline(9012,4104)(9012,4554)
\drawline(8112,4554)(8562,4554)
\drawline(8442.000,4524.000)(8562.000,4554.000)(8442.000,4584.000)
\drawline(8562,4554)(9012,4554)
\drawline(9012,5454)(9012,5004)
\drawline(8982.000,5124.000)(9012.000,5004.000)(9042.000,5124.000)
\drawline(9012,5004)(9012,4554)
\drawline(9582.000,4584.000)(9462.000,4554.000)(9582.000,4524.000)
\drawline(9462,4554)(9912,4554)
\drawline(9462,4554)(9012,4554)
\drawline(942.000,4884.000)(912.000,5004.000)(882.000,4884.000)
\drawline(912,5004)(912,4554)
\drawline(912,5004)(912,5454)
\drawline(942.000,3984.000)(912.000,4104.000)(882.000,3984.000)
\drawline(912,4104)(912,3654)
\drawline(912,4104)(912,4554)
\drawline(12,4554)(462,4554)
\drawline(342.000,4524.000)(462.000,4554.000)(342.000,4584.000)
\drawline(462,4554)(912,4554)
\drawline(912,4554)(1362,4554)
\drawline(1242.000,4524.000)(1362.000,4554.000)(1242.000,4584.000)
\drawline(1362,4554)(1812,4554)
\drawline(2712,4554)(3162,4554)
\drawline(3042.000,4524.000)(3162.000,4554.000)(3042.000,4584.000)
\drawline(3162,4554)(3612,4554)
\drawline(3612,4554)(4062,4554)
\drawline(3942.000,4524.000)(4062.000,4554.000)(3942.000,4584.000)
\drawline(4062,4554)(4512,4554)
\drawline(3612,4554)(3612,4104)
\drawline(3582.000,4224.000)(3612.000,4104.000)(3642.000,4224.000)
\drawline(3612,4104)(3612,3654)
\drawline(3612,5454)(3612,5004)
\drawline(3582.000,5124.000)(3612.000,5004.000)(3642.000,5124.000)
\drawline(3612,5004)(3612,4554)
\drawline(6342.000,4884.000)(6312.000,5004.000)(6282.000,4884.000)
\drawline(6312,5004)(6312,4554)
\drawline(6312,5004)(6312,5454)
\drawline(6882.000,4584.000)(6762.000,4554.000)(6882.000,4524.000)
\drawline(6762,4554)(7212,4554)
\drawline(6762,4554)(6312,4554)
\drawline(6312,1404)(5862,1404)
\drawline(5982.000,1434.000)(5862.000,1404.000)(5982.000,1374.000)
\drawline(5862,1404)(5412,1404)
\drawline(6312,504)(6312,954)
\drawline(6342.000,834.000)(6312.000,954.000)(6282.000,834.000)
\drawline(6312,954)(6312,1404)
\drawline(4512,1404)(4062,1404)
\drawline(4182.000,1434.000)(4062.000,1404.000)(4182.000,1374.000)
\drawline(4062,1404)(3612,1404)
\drawline(3612,1404)(3162,1404)
\drawline(3282.000,1434.000)(3162.000,1404.000)(3282.000,1374.000)
\drawline(3162,1404)(2712,1404)
\drawline(3612,1404)(3612,1854)
\drawline(3642.000,1734.000)(3612.000,1854.000)(3582.000,1734.000)
\drawline(3612,1854)(3612,2304)
\drawline(3612,504)(3612,954)
\drawline(3642.000,834.000)(3612.000,954.000)(3582.000,834.000)
\drawline(3612,954)(3612,1404)
\drawline(882.000,1074.000)(912.000,954.000)(942.000,1074.000)
\drawline(912,954)(912,1404)
\drawline(912,954)(912,504)
\drawline(882.000,1974.000)(912.000,1854.000)(942.000,1974.000)
\drawline(912,1854)(912,2304)
\drawline(912,1854)(912,1404)
\drawline(1812,1404)(1362,1404)
\drawline(1482.000,1434.000)(1362.000,1404.000)(1482.000,1374.000)
\drawline(1362,1404)(912,1404)
\drawline(912,1404)(462,1404)
\drawline(582.000,1434.000)(462.000,1404.000)(582.000,1374.000)
\drawline(462,1404)(12,1404)
\drawline(8982.000,1074.000)(9012.000,954.000)(9042.000,1074.000)
\drawline(9012,954)(9012,1404)
\drawline(9012,954)(9012,504)
\drawline(9012,1404)(8562,1404)
\drawline(8682.000,1434.000)(8562.000,1404.000)(8682.000,1374.000)
\drawline(8562,1404)(8112,1404)
\drawline(9042.000,1734.000)(9012.000,1854.000)(8982.000,1734.000)
\drawline(9012,1854)(9012,1404)
\drawline(9012,1854)(9012,2304)
\drawline(9012,1404)(9462,1404)
\drawline(9342.000,1374.000)(9462.000,1404.000)(9342.000,1434.000)
\drawline(9462,1404)(9912,1404)
\drawline(6312,1404)(6762,1404)
\drawline(6642.000,1374.000)(6762.000,1404.000)(6642.000,1434.000)
\drawline(6762,1404)(7212,1404)
\drawline(6282.000,1974.000)(6312.000,1854.000)(6342.000,1974.000)
\drawline(6312,1854)(6312,2304)
\drawline(6312,1854)(6312,1404)
\put(3162,3204){\makebox(0,0)[lb]{\smash{{{\SetFigFont{12}{14.4}{\rmdefault}{\mddefault}{\updefault}b1(u)}}}}}
\put(5862,3204){\makebox(0,0)[lb]{\smash{{{\SetFigFont{12}{14.4}{\rmdefault}{\mddefault}{\updefault}c1(u)}}}}}
\put(8562,3204){\makebox(0,0)[lb]{\smash{{{\SetFigFont{12}{14.4}{\rmdefault}{\mddefault}{\updefault}d1(u)}}}}}
\put(462,54){\makebox(0,0)[lb]{\smash{{{\SetFigFont{12}{14.4}{\rmdefault}{\mddefault}{\updefault}a2(u)}}}}}
\put(3162,54){\makebox(0,0)[lb]{\smash{{{\SetFigFont{12}{14.4}{\rmdefault}{\mddefault}{\updefault}b2(u)}}}}}
\put(5862,54){\makebox(0,0)[lb]{\smash{{{\SetFigFont{12}{14.4}{\rmdefault}{\mddefault}{\updefault}c2(u)}}}}}
\put(8562,54){\makebox(0,0)[lb]{\smash{{{\SetFigFont{12}{14.4}{\rmdefault}{\mddefault}{\updefault}d2(u)}}}}}
\put(462,3204){\makebox(0,0)[lb]{\smash{{{\SetFigFont{12}{14.4}{\rmdefault}{\mddefault}{\updefault}a1(u)}}}}}
\end{picture}
}

\end{figure}

We define the row-to-row transfer matrix on a closed chain
as $t(u) = Tr_0 \cT(u)$, where
$\cT(u)$ is the monodromy matrix given by
\begin{equation}
 \cT(u) = \cR_{01}(u)\cR_{02}(u)\cdots \cR_{0L}(u) \;.
 \label{eq:Monod}
\end{equation}
The Yang--Baxter algebra satisfied by $\cR$ ensures that transfer
matrices with different spectral parameters commute, i.e.
\begin{equation}
 [Tr_0\, \cT(u), Tr_0\, \cT(v)]=0\ , \qquad \forall u,v \;.
 \label{eq:commT}
\end{equation}
This commutativity relies on the so-called ``rail-way'' proof.

An integrable vertex model with open boundary conditions can also
be defined using a double-row monodromy matrix, cf. \cite{Skl,MN},
(see also \cite{FaWa,YaSa,KuLS,YZG})
\begin{equation}
 \cT(u) = \cR_{01}(u)\cR_{02}(u)\cdots \cR_{0L}(u)
 K(u) \cR_{L0}(u)\cR_{L-1,0}(u)\cdots \cR_{10}(u) \;.
 \label{eq:Monoddle}
\end{equation}
and transfer matrix $t(u) = Tr_0 K'(u) \cT(u)$. where $K$ and $K'$
are boundary reflection matrices. One should be able to prove that
the double-row transfer matrices commute among themselves for any
values of the spectral parameters, at least in the case $K(u)=1$
and $K'(u)=1$. This would use the reflection equation \cite{Che}
\begin{equation}
 \cR_{12}(u-v) \cK_1^-(u) \cR_{21}(u+v) \cK_2^-(v)
 =
 \cK_2^-(v) \cR_{12}(u+v) \cK_1^-(u) \cR_{21}(u-v)
 \label{eq:RE-}
\end{equation}
and a so-called crossing-unitarity relation for $R(u)$.

The commutativity for open chain has been checked with the computer
for some values of $L$ (up to now $L=2,\cdots,6$.)

In this exotic model
the weights cannot be all non-negative except for the  trivial limit u=1. So
negative weights have to be  suitably interpreted, e.g., as in \cite{Pes}.
We leave this for future investigations.

\subsection{An integrable spin chain}
\label{sect:chain}

\subsubsection{The model}
A Hamiltonian of spin chain can be defined as the first term in the
expansion of the closed transfer matrix around $u=1$
\begin{equation}
 \cH_{\rm per} =
 \left.\frac{d}{du}
 \right|_{u=1} \cT(u) . \cT(1)^{-1}
 = \sum_{i=1}^{L-1} \cH_{i\ i+1} + \cH_{L\, 1} \;,
 \label{eq:Hinteg}
\end{equation}
(for closed boundary conditions).

This Hamiltonian by construction commutes with the transfer matrices
$\cT(u)$ for any $u$.
It has a high degeneracy of spectrum (experimental-computer fact).

Similarly, an open chain Hamiltonian can be defined using the
derivative of the open chain transfer matrix
\begin{equation}
 \cH_{\rm open} =
 \left.\frac{d}{du}
 \right|_{u=1} \cT(u) . \cT(1)^{-1}
 = \sum_{i=1}^{L-1} \cH_{i\ i+1} \;,
 \label{eq:Hopen}
\end{equation}
(up to normalisation).

Let us recall the elements $B,C,D$ of the standard dual $s03$ of
$S03$ from \cite{ACDM2}. They are duals to the elements $b,c,d$,
resp., (cf. \eqref{rtt}), and their two-dimensional representation
is related to the standard sigma matrices:
\begin{equation}
 \label{eq:BCD}
 B=
 \left(
 \begin{array}{cc}
 0 & 1 \cr
 1 & 0
 \end{array}
 \right) = \sigma_1 \ ,
 \qquad
 C=
 \left(
 \begin{array}{cc}
 0 & 1 \cr
 -1 & 0
 \end{array}
 \right) = -i \sigma_2 \ ,
 \qquad
 D=
 \left(
 \begin{array}{cc}
 1 & 0 \cr
 0 & -1
 \end{array}
 \right) = \sigma_3 \  .
\end{equation}
Then the Hamiltonian on two sites may be written as: $H=B \otimes C$.
The operators $B$, $C$, $D$ satisfy
\begin{equation}
 \label{eq:BCDrel}
 B^2 = -C^2 = D^2 = 1, \quad
 DB = -BD = C, \quad
 DC = -CD = B, \quad
 CB = -BC = D
\end{equation}

\subsubsection{Eigenstates and eigenvalues}
The eigenvalues $\lambda$ and eigenstates of the Hamiltonian on two
sites are
\begin{equation}
 \label{eq:eig2states}
 \begin{array}{ll}
 \lambda =-i \qquad &
 \left| \up \up \right> + i \left| \dw\dw\right> \\[2mm]
 &
 \state{\dw \up} + i \state{\up\dw} \\[2mm]
 \lambda =i \qquad &
 \left| \up \up \right> - i \left| \dw\dw\right> \\[2mm]
 &
 \state{\dw \up} - i \state{\up\dw}
 \end{array}
\end{equation}

The eigenvalues $\lambda$ and eigenstates of the Hamiltonian on three
sites are
\begin{equation}
 \label{eq:eig3states}
 \begin{array}{ll}
 \lambda =\pm i\sqrt{2} \qquad &
 w_1^\pm=\state{\dw\dw\up} + \state{\up\dw\dw} \pm i\sqrt{2}\state{\up\up\up}
 \\[2mm]
 &
 w_2^\pm=\state{\dw\dw\up} - \state{\up\dw\dw} \mp i\sqrt{2}\state{\dw\up\dw}
 \\[2mm]
 &
 w_3^\pm=\state{\up\up\dw} + \state{\dw\up\up} \mp i\sqrt{2}\state{\dw\dw\dw}
 \\[2mm]
 &
 w_4^\pm=\state{\up\up\dw} - \state{\dw\up\up} \pm i\sqrt{2}\state{\up\dw\up}
 \end{array}
\end{equation}

The characteristic polynomial of the open chain Hamiltonian is
\begin{equation}
 \label{eq:charpoly}
 \begin{array}{ll}
 2\mbox{ sites~:}\qquad &
 (x^2+1)^2
 \\[2mm]
 3\mbox{ sites~:}\qquad &
 (x^2+2)^4
 \\[2mm]
 4\mbox{ sites~:}\qquad &
 (x^2+5)^4 (x^2+1)^4
 \\[2mm]
 5\mbox{ sites~:}\qquad &
 (x^4 + 8x^2+4)^8
 \\[2mm]
 6\mbox{ sites~:}\qquad &
 (x^2+1)^8 (x^6 + 19x^4 + 83x^2+1)^8
 \end{array}
\end{equation}

The characteristic polynomial of the periodic chain Hamiltonian is
\begin{equation}
 \label{eq:charpolyp}
 \begin{array}{ll}
 2\mbox{ sites~:}\qquad &
 x^2 (x^2+4)
 \\[2mm]
 3\mbox{ sites~:}\qquad &
 (x^2+3)^4
 \\[2mm]
 4\mbox{ sites~:}\qquad &
 x^4 (x^2+8)^2 (x^2+4)^4
 \\[2mm]
 5\mbox{ sites~:}\qquad &
 (x^4 + 10x^2+5)^8
 \\[2mm]
 6\mbox{ sites~:}\qquad &
 x^{24} (x^2+16)^4 (x^2+4)^8 (x^2+12)^8
 \end{array}
\end{equation}

\subsubsection{Symmetries}
It can be checked that the Hamiltonian on two sites commutes with
\begin{equation}
 \label{eq:symm2sites}
 B\otimes 1, \qquad
 D\otimes B, \qquad
 D\otimes D
\end{equation}
and the algebra generated by those, including
$1\otimes C$, $C\otimes D$, $B\otimes C$, $C\otimes B$.

The Hamiltonian on $L$ sites (with open boundary conditions) hence
commutes with
\begin{equation}
 \label{eq:symmL}
 B_i \equiv D \otimes D\otimes \cdots D\otimes B \otimes 1 \otimes
 \cdots \otimes 1
\end{equation}
with $B$ on $i$-th position
and with
\begin{equation}
 \label{eq:symmLa}
 B_{L+1} \equiv D \otimes D\otimes \cdots D\otimes D \otimes D \otimes
 \cdots \otimes D
\end{equation}

These $L+1$ operators generate a \emph{Clifford} algebra $\cC_{L+1}$,
i.e.
\begin{equation}
 \label{eq:Cliff}
 \{B_i, B_j\} = 2 \delta_{ij} \qquad i,j\in\{1,...,L+1\}\ .
\end{equation}

For even $L$, we have a Casimir given by
\begin{equation}
 \label{eq:casimir}
 \cC = \prod_{j=1}^{L+1} B_j \;.
\end{equation}

The left regular representation of $\cC_{L+1}$, of dimension
$2^{L+1}$ and with basis elements $B_1^{n_1} B_2^{n_2}\cdots
B_{L+1}^{n_{L+1}}$ ($n_j=0,1$) decomposes into $2^{[\frac{L+2}{2}]}$
irreducible representations, each of dimension $2^{[\frac{L+1}{2}]}$
.

For \emph{even} $L$, an irreducible representation can be described
with the following set of basis vectors
\begin{equation}
 \label{eq:basisvect}
 B_1^{n_1} B_2^{n_2} \cdots B_{[\frac{L}{2}]}^{n_{[\frac{L}{2}]}}
 \left(1 + \alpha B_{L+1}\right)
 \left(1+i \epsilon_1 B_1 B_L\right) \left(1+i \epsilon_2 B_2 B_{L-1}\right)
 \cdots
 \left(1+i \epsilon_{\frac{L}{2}} B_{\frac{L}{2}} B_{\frac{L}{2}+1}\right)
\end{equation}
on which the $B_j$ act by left multiplication.
The parameters $\epsilon_j$ satisfy $\epsilon_j^2=1$.
The exponents $n_j$ take the values~0 and~1.
On this representation, the Casimir operator $\cC$ acts as
$\prod_{j=1}^{L/2} (-i \epsilon_j)$.

We can use properties like
\begin{equation}
 \label{eq:prop}
 \left(1+i \epsilon_1 B_1 B_L\right) \left(1+i \epsilon_2 B_2 B_{L-1}\right)
 = \left(1+i \epsilon_1 B_1 B_L\right)
 \left(1- \epsilon_1 \epsilon_2 B_1 B_L B_2 B_{L-1}\right)
\end{equation}
to change the expressions, in particular to get an explicit appearance
of
$$
1 - \prod_{j=1}^{L/2} (-i \epsilon_j) B_1 B_2 \cdots B_{L+1}
\ = \ 1 - \prod_{j=1}^{L/2} (-i \epsilon_j) \cC
$$

For \emph{odd} $L$, an irreducible representation can be described
with the following set of basis vectors
\begin{eqnarray}
 \label{eq:basisvecto}
 \left(1 + (-1)^{n_{L+1}} B_{L+1}\right)
 B_1^{n_1} B_2^{n_2} \cdots B_{[\frac{L-1}{2}]}^{n_{[\frac{L-1}{2}]}}
 \left(1+i \epsilon_1 B_1 B_L\right) \left(1+i \epsilon_2 B_2 B_{L-1}\right)
 \cdots\nonumber\\
 \cdots
 \left(1+i\epsilon_{\frac{L-1}{2}} B_{\frac{L-1}{2}} B_{\frac{L+3}{2}}\right)
 \left(1+i \epsilon_{\frac{L+1}{2}} B_{\frac{L+1}{2}} B_{L+1}\right)
\end{eqnarray}

For even L denote the basis vector:

\begin{eqnarray} A^{n_1,...,n_{L\over2}} &=&
\left\{\Pi^{L\over2}_{k=1} B^{n_k}_k \right\} (1+\alpha B_{L+1})\
\Pi^{L\over2}_{j=1} (1+i\epsilon_j B_j
B_{L+1-j}) \end{eqnarray}
where $\alpha^2 =1$.

Then for $j=1,...,{L\over2}$ we have:
\begin{eqnarray} B_j A^{n_1,..,n_j,..,n_{L\over2}} &=& (-1)^{n_1+
...+n_{j-1}} A^{n_1,..,n_j+1,..,n_{L\over2}} \ ,\nn\\
 B_{L+1-j} A^{n_1,..,n_j,..,n_{L\over2}} &=&
-i\epsilon_j (-1)^{n_1+...+n_j}A^{n_1,..,n_j+1,..,n_{L\over2}}
\ ,\nn\\
 B_{L+1} A^{n_1,..,n_j,..,n_{L\over2}} &=& \alpha
(-1)^{n_1+ ...+n_{L\over2}} A^{n_1,..,n_j,..,n_{L\over2}} \ . \end{eqnarray}
The action of the Casimir is
\begin{equation} \cC\, A^{n_1,..,n_j,..,n_{L\over2}} = \alpha
\Pi^{L\over}_{j=1} (-i\epsilon_j) A^{n_1,..,n_j,..,n_{L\over2}} \end{equation}

For odd L we denote the basis vector:
\begin{eqnarray} A^{n_{L+1};n_1, .. ,n_j, .. n_{{L-1}\over 2}} &=&
(1+(-1)^{n_{L+1}}B_{L+1})\Pi^{{L-1}\over2}_{k=1}B^{n_k}_k \times
\\ && \times \left\{ \Pi^{{L-1}\over2}_{j=1}(1+i\epsilon_j B_j
B_{L+1-j})\right\}
(1+i\epsilon_{{L+1}\over2}B_{{L+1}\over2}B_{L+1}) \nonumber \end{eqnarray}

Then for $j=1,...,{{L-1}\over 2}$ we have:
\begin{eqnarray} B_j A^{n_{L+1};n_1, .. ,n_j, .. n_{{L-1}\over 2}}
&=& (-1)^{n_1 + ... n_{j-1}}A^{n_{L+1}+1;n_1, .. ,n_j+1, .. n_{{L-1}\over 2}}
\ ,\nn\\
  B_{L+1-j}A^{n_{L+1};n_1, .. ,n_j, ..
n_{{L-1}\over 2}} &=& -i\epsilon_j (-1)^{n_1 + ...
n_j}A^{n_{L+1}+1;n_1, .. ,n_j+1, .. n_{{L-1}\over 2}}\ ,\nn\\
 B_{{L+1}\over2}A^{n_{L+1};n_1, .. ,n_j, ..
n_{{L-1}\over 2}} &=& i\epsilon_{{L+1}\over2} (-1)^{n_{L+1}+1}
A^{n_{L+1}+1;n_1, .. ,n_j, .. n_{{L-1}\over 2}}\ ,\nn\\
 B_{L+1}A^{n_{L+1};n_1, .. ,n_j, ..
n_{{L-1}\over 2}} &=& (-1)^{n_{L+1}}A^{n_{L+1};n_1, .. ,n_j, ..
n_{{L-1}\over 2}} \ . \end{eqnarray}

In the derivation of the above relations were used also the
following formulae:

For even L:
\begin{eqnarray} B_{L+1-j}(1 +i\epsilon_j B_j B_{L+1-j}) &=&
-i\epsilon_j B_j (1 +i\epsilon_j B_j B_{L+1-j})\ , \nn\\
 B_{L+1}(1+\alpha B_{L+1}) &=& \alpha (1+\alpha
 B_{L+1}) \ . \end{eqnarray}

For odd L:
\begin{eqnarray} B_{L+1-j} (1+i\epsilon_j B_j B_{L+1-j}) &=& -i
\epsilon_j B_j (1+i\epsilon_j B_j B_{L+1-j})\ , \nn\\
 B_{{L+1}\over2} (1+i\epsilon_{{L+1}\over2}
B_{{L+1}\over2} B_{L+1}) &=& i\epsilon_{{L+1}\over2} B_{L+1}
(1+i\epsilon_{{L+1}\over2} B_{{L+1}\over2} B_{L+1})\ , \nn\\
 B_{L+1} (1+(-1)^{n_{L+1}} B_{L+1}) &=&
(-1)^{n_{L+1}} (1+(-1)^{n_{L+1}} B_{L+1})\ . \end{eqnarray}

\section*{Acknowledgments}
The work of DA, VKD and SGM was supported in part by
 the TMR Network EUCLID, contract HPRN-CT-2002-00325.
The work of VKD and SGM was supported in part also by the
Bulgarian National Council for Scientific Research, grant
F-1205/02,   and the European RTN  'Forces-Universe', contract
MRTN-CT-2004-005104. We would like also to thank the anonymous
referee for useful remarks.


\end{document}